\documentclass{amsproc}

\usepackage{mathtools,amssymb,url}
\usepackage[boxsize=1.2em]{ytableau}

\usepackage{mathrsfs} 

\numberwithin{equation}{section}

\newcommand{\bq}[1]{\boldsymbol{(#1)}}
\renewcommand{\a}{\alpha}
\renewcommand{\b}{\beta}
\renewcommand{\c}{\gamma}
\DeclareMathOperator{\PF}{PF}
\DeclareMathOperator{\des}{des}
\DeclareMathOperator{\maj}{maj}
\newcommand{\qbinom}[2]{\genfrac{[}{]}{0pt}{}{#1}{#2}}
\newcommand{\N}{\mathbb{N}}
\newcommand{\Z}{\mathbb{Z}}
\newcommand{\R}{\mathbb{R}}
\renewcommand{\P}{\mathbb{P}}
\renewcommand{\o}{\omega}
\newcommand{\s}{\sigma}
\newcommand{\A}{\mathscr{A}}

\renewcommand{\L}{\mathscr{L}}
\renewcommand{\S}{\mathscr{S}}
\newcommand{\sP}{\mathscr{P}}

\begin{document}

\title[$P$-Partitions]{A Historical Survey of  $P$-Partitions}

\author{Ira M. Gessel}
\address{Brandeis University}
\email{gessel@brandeis.edu}

\subjclass[2010]{Primary 05A15. Secondary 05A17}

\date{June 2, 2015}

\begin{abstract}
We give a historical survey of the theory $P$-partitions, starting with MacMahon's work, describing Richard Stanley's contributions and his related work, and continuing with more recent developments.
\end{abstract}

\maketitle


\emph{Dedicated to Richard Stanley on the occasion of his seventieth birthday.}
\section{Introduction}
Richard Stanley's 1971 Harvard Ph.D.~thesis \cite{stanley-thesis}  studied two related topics, plane partitions and $P$-partitions. A short article on his work on $P$-partitions appeared in  1970 
\cite{chromatic}, and a detailed exposition of this thesis work appeared as an American Mathematical Society Memoir \cite{ordered} in 1972. In this paper, I will describe some of the historical background of the theory of $P$-partitions, sketch Stanley's contribution to the theory, and mention some more recent developments.

The basic idea of the theory of $P$-partitions was discovered by P. A. MacMahon, generalized by Knuth, and independently rediscovered by Kreweras. In order to understand the different notations and approaches that these authors used, it will be helpful to start with a short exposition of Stanley's approach.

\section{Stanley's theory of $P$-partitions}
\label{s-theory}

I will usually follow Stanley's notation, but with some minor modifications. Stanley has given a more recent account of the theory of $P$-partitions in \emph{Enumerative Combinatorics}, Vol.~1 \cite{ec1}, Sections 1.4 and 3.15.

\begin{figure}[thbp] 
   \centering
   \includegraphics[width=1in]{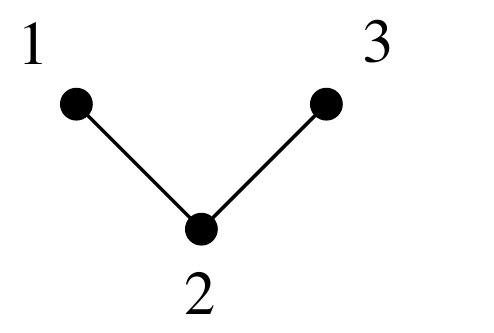} 
   \caption{A poset with three elements}
   \label{f-1}
\end{figure}

Let $P$ be a finite partially ordered set (poset) with $p$ elements and let $\o$ be a bijection from $P$ to $[p]=\{1,2,\dots, p\}$, called a \emph{labeling} of $P$. We use the symbols $\preceq$ and $\prec$ for the partial order relation of $P$. Then a $(P,\o)$-partition  is a map $\sigma$ from $P$ to the set $\N$ of nonnegative integers satisfying the conditions
\begin{enumerate}
\item[(i)] If $X\prec Y$  then $\s(X)\ge\s(Y)$.
\item[(ii)] If $X\prec Y$ and $\o(X)>\o(Y)$  then $\s(X)> \s(Y)$.
\end{enumerate}
If $\sum_{X\in P}\s(X) = n$ then we call $\s$ a $P$-partition of $n$.
For simplicity, we assume that the elements of $P$ are $1,2,\dots, p$.
We denote by $\A(P,\o)$ the set of all $(P,\o)$-partitions and by $\A(P,\o;m)$ the set of all $(P,\o)$-partitions with largest part at most $m$.
Thus for the labeled poset $P$ of Figure \ref{f-1}, where the elements of $P$ are identified with   their labels, $\sigma: \{1,2,3\} \to \N$ is a $P$-partition if and only if $\s(2)>\s(1)$ and $\s(2)\ge\s(3)$.

If $\o(X)<\o(Y)$ whenever $X\prec Y$, then $\o$ is called a \emph{natural labeling}, and if 
 $\o(X)>\o(Y)$ whenever $X\prec Y$, then $\o$ is called a \emph{strict labeling}.
 If $\o$ is a natural labeling, then all of the inequalities in the definition of a $P$-partition are weak, and if $\o$ is a strict labeling then all the inequalities are strict.

A \emph{linear extension} of $P$ is a total order (or chain) on $P$ that contains $P$. Every linear extension of $P$ inherits the labeling $\o$ of $P$. We denote by $\L(P,\o)$ the set of linear extensions of $P$. We may identify a linear extension of $P$  with the  permutation of $[p]$ obtained by listing the labels of its elements in increasing order. Thus for $P$ in Figure \ref{f-1}, $\L(P,\o)$ consists of the two permutations $213$ and $231$. Then the ``fundamental theorem of $P$-partitions'' (stated somewhat differently, but equivalently, by Stanley \cite[Lemma 6.1]{ordered}) asserts that the set $\A(P,\o)$ of $P$-partitions is the disjoint union of $\A(\pi,\omega)$ over all elements $\pi$ of $\L(P,\o)$.

Thus for the poset $P$ of Figure \ref{f-1}, the fundamental theorem says that the set of $P$-partitions is the disjoint union of the solutions of $\s(2)>\s(1)\ge \s(3)$ and the solutions of $\s(2)\ge\s(3)>\s(1)$.

We sketch here Stanley's proof. Given a $P$-partition $\s$, we can arrange its values in weakly decreasing order and thus find a permutation $\pi$ of $P$ such that 
\begin{equation*}
\s(\pi(1)) \ge \s(\pi(2)) \ge \dots \ge \s(\pi(p)).
\end{equation*}
There may be many ways to do this, but we get a unique permutation if we require that the labels of equal values be arranged in increasing order; i.e., if $\s(\pi(i)) = \s(\pi({i+1}))$ then $\o(\pi(i)) < \o(\pi({i+1}))$. But the contrapositive of this  property is (since the labels are all distinct and $\s(\pi(i)) \ge \s(\pi({i+1}))$) that if $\o(\pi(i)) > \o(\pi({i+1}))$ then $\s(\pi(i)) > \s(\pi({i+1}))$. This shows that every $P$-partition is in $\A(\pi,\o)$ for a unique $\pi\in \L(P,\o)$, and it is clear from the definitions that $\A(\pi,\o)\subseteq \A(P,\o)$ for each $\pi\in \L(P,\o)$.

The fundamental theorem has enumerative consequences, since $(P,\o)$-partitions are easy to count when $P$ is a total order.   In particular, Stanley defines
\begin{align*}
U_m(P,\o) &= \sum_{\s\in \A(P,\o;m)}q^{\s(1)+\cdots+\s(p)},\\
U(P,\o) &= \lim_{m\to\infty}U_m(P,\o) =\sum_{\s\in \A(P,\o)}q^{\s(1)+\cdots+\s(p)},\\
\Omega(P,\o;m)&=U_{m-1}(P,\o)|_{q=1};
\end{align*}
$\Omega(P,\o;m)$ is called the \emph{order polynomial} of $(P,\o)$.
(In section \ref{s-recip} we will write $U_m(P,\o;q)$ and $U(P,\o;q)$ when we need to show dependence on $q$.)
To evaluate these sums we apply the fundamental theorem, which reduces the problem to the case in which $P$ is a total order.  
If $\pi$ is a permutation of $[p]$, we define the \emph{descent set} $\S(\pi)$ to be the set $\{\,j \mid \pi(j) > \pi(j+1)\,\}$, and we denote by $\des(\pi)$ the number of elements of $\S(\pi)$ and by $\maj(\pi)$ (the \emph{major index}\footnote{Stanley used the term ``index".} of $\pi$) the sum of the elements of $\S(\pi)$.
Then the basic result on these quantities is 
\begin{equation}
\label{e-majdes}
\sum_{m=0}^\infty U_m(P,\o) t^m = \frac{\sum_{\pi\in \L(P,\o)}q^{\maj(\pi)} t^{\des(\pi)}}{(1-t)(1-tq)\cdots(1-tq^p)}.
\end{equation}
Equation \eqref{e-majdes} is easily proved directly if $P$ is a total order (chain), and the general case follows from the fundamental theorem. As consequences of \eqref{e-majdes} we have
\begin{gather}
\label{e-U}
U(p,\o) = \frac{\sum_{\pi\in \L(P,\o)}q^{\maj(\pi)}}{(1-q)\cdots(1-q^p)}
\intertext{and}
\label{e-order}
\sum_{m=0}^\infty \Omega(P,\o;m) t^m = \frac{\sum_{\pi\in \L(P,\o)} t^{1+\des(\pi)}}{(1-t)^{p+1}}.
\end{gather}

\section{MacMahon}

The story of $P$-partitions begins with Percy A. MacMahon's work on plane partitions \cite{macmahon1911memoir} (see also \cite[Vol.~2, Section X, Chapter 3]{MR0141605}). The problem that MacMahon considers is that of counting plane partitions of a given shape; that is,  arrangements of nonnegative integers with a given sum in a ``lattice" such as 
\begin{equation*}
\begin{array}{cccc}
4&4&2&1\\
4&3&2\\
2&1
\end{array}
\end{equation*}
in which the entries are weakly decreasing in each row and column. 
MacMahon gives a simple example to illustrate his idea. We want to count arrays of nonnegative integers
\begin{equation*}
\begin{array}{cc}
p&q\\
r&s
\end{array}
\end{equation*}
satisfying $p\ge q\ge  s$ and $p\ge r\ge s$, and we assign to such an array the weight $x^{p+q+r+s}$.
The set of solutions of these inequalities is the disjoint union of the solution sets of the inequalities
\begin{equation*}
\textrm{(i) } p\ge q\ge r\ge s \quad \textrm{and}\quad  \textrm{(ii) } p\ge r> q\ge s.
\end{equation*}
Setting $r=s+A$, $q=s+A+B$, and $p=s+A+B+C$, where $A$, $B$, and $C$ are arbitrary nonnegative integers, we see that the sum $\sum x^{p+q+r+s}$ over solutions of 
$p\ge q\ge r\ge s $
is equal to
\begin{equation*}
\sum_{A,B,C,s\ge0}x^{C+2B+3A+4s}=\frac{1}{\bq 1\bq 2\bq 3 \bq 4},
\end{equation*}
where $\bq n=  (1-x^n)$. Similarly, the generating function for solutions of $ p\ge r> q\ge s$ is $x^2/\bq 1 \bq 2 \bq 3 \bq 4$, so the generating function for all of the arrays is 
\begin{equation}
\label{e-gf1}
\frac{1+x^2}{\bq 1\bq2 \bq3 \bq4}.
\end{equation}
MacMahon explains (but does not prove) that a similar decomposition exists for counting plane partitions of any shape, and moreover, the terms that appear in the numerator have combinatorial interpretations. They correspond to what MacMahon calls \emph{lattice arrangements}, which are essentially what we now call \emph{standard Young tableaux}. In the example under discussion there are two lattice arrangements, 
\begin{equation*}
\begin{array}{cc}
4&3\\
2&1
\end{array}
\text{ and } 
\begin{array}{cc}
4&2\\
3&1
\end{array}.
\end{equation*}
They are the plane partitions of the shape under consideration in which the entries are $1, 2, \dots, n$, where $n$ is the number of positions in the shape. To each lattice arrangement MacMahon associates a \emph{lattice permutation}:  the $i$th entry in the lattice permutation corresponding to an arrangement is the row of the arrangement in which $n+1-i$ appears, where the rows are represented by the Greek letters $\a$, $\b$, \dots.  So the lattice permutation associated to the first arrangement is $\a\a\b\b$ and to the second is $\a\b\a\b$. (A sequence of Greek letters is called a lattice permutation if any initial segment contains at least as many $\a$s as $\b$s, at least as many $\b$s as $\gamma$s, and so on.) To each lattice permutation, MacMahon associates an inequality relating $p$, $q$, $r$, and $s$; the $\a$s are replaced, in left-to-right order with the first-row variables $p$ and $q$, and the $\b$s are replaced with the second-row variables $r$ and~$s$. A greater than or equals sign is inserted between two Greek letters that are in alphabetical order and a greater than sign is inserted between two Greek letters that are out of  alphabetical order. So the lattice permutation $\a\a\b\b$ gives the inequalities $p\ge q\ge r\ge s$ and the lattice permutation $\a\b\a\b$ give the inequalities $p\ge r>q\ge s$. Each lattice permutation contributes one term to the numerator of \eqref{e-gf1}, and the power of $x$ in such a term is the sum of the positions of the Greek letters that are followed by a smaller Greek letter. MacMahon then describes the variation in which a restriction on part sizes is imposed. The decomposition into disjoint inequalities works exactly as in the unrestricted case, and reduces the problem to counting partitions with a given number of parts and a bound on the largest part. Most of the rest of the paper is devoted to applications of this idea to the enumeration of plane partitions. In a postscript, MacMahon considers the analogous situation in which only decreases in the rows are required, not in the columns. The enumeration of such arrays is not of much interest in itself, since the generating function for an array with $p_1, p_2,\dots, p_n$ nodes in its $n$ rows is clearly
\begin{equation*}
\frac{1}{\bq 1 \cdots \bq {p_1}\bq 1 \cdots \bq {p_2}\cdots\cdots \bq 1\cdots \bq {p_n}}.
\end{equation*}
However the same decomposition that is used in the case of plane partitions yields interesting results about permutations.

In a follow-up paper \cite{macmahon1913indices} (see also \cite[Vol. 2, Section IX, Chapter 3]{MR0141605}), MacMahon elaborates on this idea. Given a sequence of elements of a totally ordered set,  MacMahon defines a \emph{major contact}\footnote{Now usually called a \emph{descent}.} to be a pair of consecutive entries in which the first is greater than the second, and he defines the \emph{greater index}\footnote{MacMahon's greater index is now usually called the \emph{major index}. The term ``major index" was introduced by Foata and Sch\"utzenberger \cite{MR519777,MR506852} in reference to MacMahon's military rank.
Curiously, MacMahon used the term ``major index" for a related concept that does not seem to have been further studied.}  to be the sum of the positions of the first elements of the major contacts. 
Thus the greater index of $\b\a\a\a\c\c\b\a\c$, where the letters are ordered alphabetically, is $1+6+7=14$. (He similarly defines the ``equal index" and ``lesser index" but these do not play much of a role in what follows.)
MacMahon's main result in the paper is that the sum $\sum x^p$, where $p$ is the greater index, over all ``permutations of the assemblage $\a^i \b^j\c^k\cdots$" is 
\begin{equation*}
\frac{\bq 1\bq 2 \cdots \bq{i+j+k+\cdots}}
{\bq 1\bq 2 \cdots \bq i \cdot \bq 1\bq 2 \cdots \bq j \cdot \bq 1\bq 2 \cdots \bq k \cdots}
\end{equation*}
As in the previous paper MacMahon illustrates with an example, but does not give a formal proof, nor even an informal explanation of why the procedure works. We consider the sum of $x^{a_1+a_2+a_3+b_1+b_2}$ over all inequalities $a_1\ge a_2\ge a_3$, $b_1\ge b_2$. We see directly that the sum is 
\begin{equation*}
\frac{1}{\bq 1 \bq 2 \bq 3\cdot \bq 1 \bq 2}.
\end{equation*}
MacMahon breaks up these inequalities just as before into subsets corresponding to all the permutations of $\a^3\b^2$; for example, to the permutation $\a\b\a\b\a$ corresponds the inequalities $a_1\ge b_1 > a_2 \ge b_2 > a_3$, where the strict inequalities correspond to the major contacts, which in this example are all of the form $\b \a$. The generating function for the solutions of these inequalities is 
\begin{equation*}
\frac{x^6}{\bq 1\bq 2\bq 3\bq 4 \bq 5};
\end{equation*}
here $6=2+4$ is the the greater index of the permutation $\a\b\a\b\a$. Summing the contributions from all ten permutations of $\a^3\b^2$ gives
\begin{equation*}
\frac{\sum x^p}{\bq 1\bq 2\bq 3\bq 4 \bq 5} = \frac{1}{\bq 1 \bq 2 \bq 3\cdot \bq 1 \bq 2}.
\end{equation*}

In his book \emph{Combinatory Analysis} \cite[Vol.~2, Section IX]{MR0141605} MacMahon discusses the analogous result when a bound is imposed on the part sizes. The sum of $x^{a_1+\cdots+a_p}$ over all solutions of $n\ge a_1\ge \cdots \ge a_p$ is  
$\bq {n+1}\cdots \bq{n+p}/\bq1\cdots \bq p$, and MacMahon derives in Art.~462 an important, though not well-known, formula that he writes as
\begin{equation}
\label{e-qSN1}
\begin{multlined}
\sum_{n=0}^\infty g^n \frac{\bq{n+1}\cdots \bq{n+p_1}\cdot \bq{n+1}\cdots \bq{n+p_2}\cdots\cdots \bq{n+1}\cdots \bq{n+p_m}} 
{\bq 1 \bq 2 \cdots \bq{p_1} \cdot\bq 1 \bq 2 \cdots \bq{p_2}\cdots\cdots \bq 1 \bq 2 \cdots \bq{p_m}}\\
 = \frac{1+g\PF_1+g^2\PF_2+\cdots+g^\nu \PF_\nu}
 {(1-g)(1-gx)(1-gx^2)\cdots(1-gx^{p_1+\cdots+p_\nu})}.
\end{multlined}
\end{equation}
Here $\PF_s$ is the generating function, by greater index, of permutations of the assemblage $\a_1^{p_1}\a_2^{p_2}\cdots \a_m^{p_m}$ with $s$ major contacts. This result is worth restating it in more modern terminology: Let 
\[A_{p_1,\dots, p_m}(t,q) = \sum_\pi t^{\des(\pi)}q^{\maj(\pi)},\] 
where the sum is over all permutations $\pi$ of the multiset $\{1^{p_1}, 2^{p_2}, \dots, m^{p_m}\}$, and if $\pi=a_1\cdots a_p$, where $p=p_1+\cdots p_m$ then $\des(\pi)$ is the number of descents of $\pi$, that is, the number of indices $i$ for which $a_i>a_{i+1}$, and $\maj(\pi)$ is the sum of the descents of $\pi$. Let $(a;q)_m$ be the $q$-rising factorial $(1-a)(1-aq)\cdots (1-aq^{n-1})$, let $(q)_n$ denote $(q;q)_n = (1-q)\cdots (1-q^n)$ and let
$\qbinom{m}{n}$ denote the $q$-binomial coefficient
\begin{equation*}
\frac{(q)_{m}}{(q)_n (q)_{m-n}}.
\end{equation*}
Then 
\begin{equation}
\label{e-qSN2}
\sum_{n=0}^\infty t^n \qbinom{n+p_1}{p_1}\qbinom{n+p_2}{p_2}\cdots \qbinom{n+p_m}{p_m} 
  =\frac {A_{p_1,\dots, p_m}(t,q)}{(t;q)_{p+1}}.
\end{equation}
Several specializations of \eqref{e-qSN2} are worth mentioning. 
For $q=1$ 
 the polynomials $A_{p_1,\dots, p_m}(t,1)$ solve Simon Newcomb's problem, the problem of counting permutations of a multiset by descents.\footnote{In \cite{macmahon1908second} MacMahon describes Simon Newcomb's problem as the equivalent problem of counting permutations of a multiset by consecutive pairs which are \emph{not} descents, though in his book \cite[volume 1, pp. 187]{MR0141605} he counts by descents.} (MacMahon had solved Simon Newcomb's problem earlier by a different method \cite{macmahon1908second}; however, he does not note here the connection with Simon Newcomb's problem.) In the case $q=1$, $p_1=\cdots=p_m=1$, the polynomials $A_{1^m}(t,1)$ are the Eulerian polynomials\footnote{The Eulerian polynomials are often defined to be our $tA_{1^m}(t,1)$, making the generating function the nicer
$\sum_{n=0}^\infty t^n n^m$.
}, satisfying
\begin{equation*}
\sum_{n=0}^\infty t^n (n+1)^m = \frac{A_{1^m}(t,1)}{(1-t)^{m+1}}.
\end{equation*}
For $p_1=\cdots=p_m=1$, \eqref{e-qSN2} becomes
\begin{equation*}
\sum_{n=0}^\infty t^n (1+q+\cdots +q^n)^m
  =\frac {A_{1^m}(t,q)}{(t;q)_{m+1}},
\end{equation*}
a result often attributed to Carlitz \cite{MR0366683} (though Carlitz stated an equivalent result much earlier \cite{MR0060538},  attributing it to John Riordan).

\section{Kreweras}
In 1967, Germain Kreweras \cite{kreweras1967} (see also \cite{MR0200180} for a briefer account) used an approach similar to MacMahon's, though stated very differently, to solve a common generalization of Simon Newcomb's problem\footnote{Kreweras's seems to be unaware of MacMahon's work on Simon Newcomb's problem and refers only to Riordan \cite[216--219]{MR0096594} as a reference.}   and what he calls ``Young's problem". In Young's problem, we are given two weakly decreasing sequences $Y=(y_1,\dots, y_h)$ and $Y'=(y_1',\dots, y_h')$ with $y_i\ge y_i'$ for each $i$ and we ask how many ``Young chains" there are from $Y'$ to $Y$, which are sequences of partitions (weakly decreasing sequences of integers) starting with $Y'$ and ending with $Y$ in which each partition is obtained from the previous one by increasing one part by 1. In modern terminology, these are standard tableaux of shape $Y/Y'$; that is, fillings of a Young diagram of shape $Y$ with the squares of a Young diagram of shape $Y'$ removed from it, with the integers $1,2,\dots, m$ (where $m$ is the total number of squares), so that the entries are increasing in every row and column. For example,  if $Y'= (2,1,0)$ and $Y=(3,2,2)$ then one of the Young chains from $Y'$ to $Y$ is $Y'=(2,1,0), (3,1,0), (3,1,1), (3,2,1), (3,2,2) = Y$. This corresponds to the skew Young tableau 
\begin{equation*}
\ytableaushort{\none\none1,\none3,24}
\end{equation*}
in which the entry $i$ occurs in row $j$ if the $i$th step in the chain is an increase by 1 in the $j$th position. 

Kreweras defines a ``return" (\emph{retour en arri\`ere}) of a Young chain to consist of three consecutive partitions $UVW$ such that the entry augmented in passing from $V$ to $W$ has an index that is strictly less than which is augmented in passing from $U$ to $V$. In terms of Young tableau, a return corresponds to an entry $i$ which is in a higher row than $i+1$. (In our example, 1 and 3 correspond to returns.) In MacMahon's approach, the returns correspond to major contacts of lattice permutations. Kreweras writes $\theta_r(Y,Y')$ for the number of Young chains from $Y'$ to $Y$ with $r$ returns.

He observes that Simon Newcomb's problem is equivalent to a special case of computing $\theta_r(Y,Y')$; thus the number of  
permutations of the multiset $\{1^3, 2, 3^2\}$ with $r$ descents is equal to the number of skew Young tableaux of shape
\begin{equation*}
\ydiagram{3+3,2+1,2}
\end{equation*}
with $r$ returns.
He then gives the solution to this problem in the form
\begin{equation*}
\frac{\sum_{r\ge0} \theta_r(Y,Y') t^r}{(1-t)^{\eta-\eta'+1}}=\sum_{r\ge0} w_r t^r.
\end{equation*}
Here $\eta$ is the sum of the entries of $Y$, $\eta'$ is the sum of the entries of $Y'$, and $w_r$ is the number of chains 
\begin{equation}
\label{e-chains}
Y' \le Z_1 \le \cdots \le Z_r \le Y;
\end{equation}
 in an earlier work \cite{kreweras1965}, Kreweras had given the formula 
\begin{equation*}
w_r = \det \left(\binom{y_i-y_j'+r}{i-j+r}\right)_{i,j=1,\dots, h},
\end{equation*}
where $Y=(y_1,\dots, y_h)$ and $Y' = (y'_1,\dots, y'_h)$. 
In the case of Simon Newcomb's problem, the determinant is upper triangular, and is therefore a product of binomial coefficients (as can also be seen directly).

Kreweras's method of proof is ultimately equivalent to MacMahon's approach, though described very differently: he associates to every chain \eqref{e-chains} a Young chain from $Y'$ to $Y$ in such a way that the contribution to $\sum_r w_rt^r$ corresponding to  a given Young chain with $r$ returns is $t^r/(1-t)^{\eta-\eta'+1}$.

In a later paper \cite{MR623036}, Kreweras studied what is in Stanley's terminology the order polynomial of a naturally labeled poset. Although published in 1981, long after Stanley's memoir \cite{ordered}, Kreweras stated that Stanley's work was unknown to him when the paper was written.

\section{Knuth}

In 1970, Donald E. Knuth \cite{MR0277401} used MacMahon's approach to study solid (i.e., three-dimensional) partitions. MacMahon had conjectured that the generating function for solid partitions was $\prod_{i=1}^\infty (1-z^i)^{-\binom{i+1}{2}}$. This conjecture had been disproved earlier \cite{MR0217029}, but Knuth wanted to compute the number $c(n)$ of solid partitions of $n$ for larger values of $n$ in an (unsuccessful) attempt to find patterns. Knuth realized that MacMahon's approach would work for arbitrary partially ordered sets, not just those corresponding to plane partitions.
Knuth takes a set $P$  partially ordered by the relation $\prec$ and well-ordered by the total order $<$, where $x\prec y$ implies $x<y$. He defines a $P$-partition of $N$ to be a function $n$ from $P$ to the set  of nonnegative integers satisfying (i) $x\prec y$ implies $n(x)\ge n(y)$,  (ii) only finitely many $x$ have $n(x)>0$, and  (iii) $\sum_{x\in P} n(x) = N$.
Knuth proves that there is a bijection from $P$-partitions of $N$ to pairs of sequences
\begin{gather*}
n_1\ge n_2\ge\cdots \ge n_m\\
x_1, x_2, \dots, x_m
\end{gather*}
where $m\ge0$, the $n_i$ are positive integers with sum $N$,  and the $x_i$ are distinct elements of $P$ satisfying
\begin{enumerate}
\item[(S1)] For $1\le j\le m$ and $x\in P$, $x\prec x_i$ implies $x=x_i$ for some $i<j$.
\item[(S2)] $x_i>x_{i+1}$ implies $n_i > n_{i+1}$ for $1\le i < m$. 
\end{enumerate}
Knuth is interested primarily in the case in which $P$ is countably infinite, for which he uses a modification of this bijection to prove that if $P$ is an infinite poset and $s(n)$ is the number of $P$-partitions of $n$ then 
\begin{equation*}
1+s(1)z+s(2)z^2 + \cdots  = (1+t(1)z+t(2)z ^2+\cdots)/(1-z)(1-z^2)(1-z^3)\cdots
\end{equation*}
where $t(k)$ is the number of linear extensions of finite order ideals of $P$ with ``index" $k$;  Knuth's index is a variant of MacMahon's greater index.

\section{Thomas}

Gl\^anffrwd Thomas's 1977 paper \cite{thomas}, based on his 1974 Ph.D. thesis \cite{thomas-thesis} appeared after Stanley's memoir, but it was written without knowledge of Stanley's work (but with knowledge of MacMahon's). Thomas's motivation was the combinatorial definition of Schur functions. If $\lambda$ is a  partition, then a  Young tableau of shape $\lambda$ is a filling of the Young diagram of $\lambda$ that is weakly increasing in rows and strictly increasing in columns. For example, if $\lambda$ is the partition $(4,2,1)$ then a Young tableau of shape $\lambda$ is 
\begin{equation}
\label{e-ssyt}
\ytableaushort{4411,21,1}
\end{equation}
The \emph{Schur function} $s_\lambda$ is the sum of the weights of all Young tableaux of shape $\lambda$, where the weight of a Young tableau is the product of $x_i$ over all of its entries $i$. (So the weight of the tableau \eqref{e-ssyt}  is $x_1^4x_2x_4^2$.) Schur functions are symmetric in the variables $x_i$ and have important applications in enumeration and in the representation theory of symmetric and general linear groups.

Thomas considers a more general situation, in which he allows as  shapes (which he calls ``frames") any subset of $\Z\times \Z$ and he defines a \emph{numbering} of a frame to be filling with positive integers that is weakly increasing in rows and strictly increasing in columns. For example, 
\begin{equation}
\label{e-frame}
\ytableaushort{12,\none4,2\none4}
\end{equation}
is a numbering. To any numbering he associates an \emph{index numbering} by replacing its entries in increasing order with $1,2,\dots, m$, where $m$ is the number of entries, and ties are broken from bottom to top and then left to right.
Thus the index numbering corresponding to \eqref{e-frame} is 
\begin{equation*}
\ytableaushort{13,\none5,2\none4}
\end{equation*}
Thomas calls two numberings equivalent if they have the same index numbering. He defines the monomial of  a numbering of a frame to be the product $\prod x_i$ over all the entries $i$ of the frame. (So the sum of the monomials of all the frames of a Young diagram is a Schur function.) His goal is determine the sum of the monomials of an equivalence class of numberings.

The numberings of frames are a particular case of $P$-partitions corresponding to subposets of $\Z\times \Z$, and except for the case of  skew Schur functions, which are symmetric, one might just as well study general $P$-partitions with his weighting. The interesting aspect of his work is that it seems to be the earliest appearance (after a brief mention by Stanley \cite[p.~81]{ordered}) of what are now called quasi-symmetric generating functions for $P$-partitions, which we will discuss in more detail in Section\ref{s-qs}.

As we have seen, the study of $P$-partitions leads to inequalities like $j_1\ge j_2 > j_3 \ge j_4$, or equivalently (following Thomas),
\[i_1\le i_2< i_3 \le i_4.\]
While MacMahon wanted to compute $\sum x^{i_1+\cdots +i_4},$  Thomas was interested in the more refined multivariable generating function
\begin{equation*}
\sum_{i_1\le i_2< i_3 \le i_4} x_{i_1}x_{i_2}x_{i_3}x_{i_4}.
\end{equation*}
This is  a \emph{fundamental quasi-symmetric function}; these form a basis for the algebra of quasi-symmetric functions, which will be discussed in Section \ref{s-qs}. 

Thomas  applies \emph{Baxter operators} to construct quasi-symmetric functions.
A Baxter operator on a commutative algebra $A$ over a field $K$ is a linear operator $B: A\to A$ such that for some fixed nonzero $\theta\in K$, 
\begin{equation*}
B(aB(b)) + B(bB(a)) = B(a)B(b) + B(\theta ab)
\end{equation*} 
for all $a,b\in A$. 
Now let $A$ be the algebra of infinite sequences $(a_1,a_2,\dots)$ with entries in a field, with componentwise operations. 

We define two maps $A\to A$; first a map introduced by Rota and Smith \cite{MR0343094}
\begin{align*}
S(a_1,a_2,\dots) &= \biggl(0,a_1, a_1+a_2,\dots, \sum_{i=1}^{r-1}a_i, \dots\biggr),
\end{align*}
which we write as
$\left(\dots, \sum_{i=1}^{r-1}a_i, \dots\right)$,
and a variant
\begin{align*}
P(a_1,a_2,\dots)&= \biggl(a_1, a_1+a_2,\dots, \sum_{i=1}^{r}a_i, \dots\biggr)\\
 &=\biggl(\dots, \sum_{i=1}^{r}a_i, \dots\biggr).
\end{align*}
Then $S$ and $P$ are Baxter operators.

We show by an example  the connection between these operators and quasi-symmetric functions:
Let $x=(x_1, x_2, x_3,\dots)$. Then 
\begin{equation*}
xS(xS(xP(x))) = \biggl(\dots, \sum_{1\le i\le j<k<r} x_i x_j x_k x_r,\dots \biggr).
\end{equation*}
Thus the fundamental quasi-symmetric function 
\begin{equation*}
\sum_{i\le j<k<l} x_i x_j x_k x_l,\dots 
\end{equation*}
is obtained by adding all the entries of $xS(xS(xP(x)))$. 

Thomas's approach has been further developed by Rudolf Winkel \cite{MR1612387}.

\section{Stanley}

In this section we discuss a few highlights of Stanley's memoir \cite{ordered}. 

\subsection{Reciprocity theorems}
\label{s-recip}
If $\o$ is a labeling of a poset $P$ of size $p$, we  define the \emph{complementary labeling} $\bar\o$ of $P$ by $\bar\o(i) = p+1-\o(i)$. Thus when we change the labeling of $P$ from $\o$ to $\bar\o$, the strict and weak inequalities in the definition of a $P$-partition are switched. If the permutation $\pi$ in $\L(P,\o)$ corresponds to the permutation $\bar\pi$ in $\L(P, \bar\o)$ then the descent sets $\S(\pi)$ and $\S(\bar\pi)$ are complementary subsets of $[p-1]$, and thus by \eqref{e-majdes}--\eqref{e-order}, $U_m(P,\bar\o;q)$, $ U(P, \bar\o;q)$, and $\Omega(P,\bar\o;m)$ are determined by $U_m(P,\o;q)$, $ U(P, \o;q)$, and $\Omega(P,\o;m)$. The formulas expressing these relations are surprisingly simple. We first note that if $P$ is a chain and $(P,\o)$ corresponds to a permutation $\pi$ with $s$ descents then 
$U_m(P, \o;q)= q^{\maj(\pi)}\qbinom{p+m-s}{p}$. It follows that in general, $U_m(P,\o;q)$ is a polynomial in $q^m$ and $\Omega(P,\o;m)$ is a polynomial in $m$,  so $U_m(P,\o;q)$ and $\Omega(P,\o; m)$ can be extended in a natural way to negative values of $m$. Moreover, $U(P,\o;q)$ is a rational function of $q$, so $U(P,\o; q^{-1})$ is well-defined as a rational function of $q$. Then we have the following reciprocity formulas
\begin{align}
\notag
q^p U_m(P,\bar\o;q) &= (-1)^p U_{-(m+2)}(P,\o, q^{-1}), \text{ with $U_{-1}(P,\o) = 0$}\\
\notag
q^p U(P, \bar\o;q) &= (-1)^q U(P, \o; q^{-1})\\
\label{e-rec-ord}
\Omega(P,\bar\o;m) &= (-1)^p\Omega(P,\o; -m).
\end{align}
When $P$ satisfies certain ``chain conditions" there is an additional relation between the enumerative quantities associated with $(P,\o)$ and $(P,\bar\o)$ \cite[18--19]{ordered}. We state here one of these results \cite[Proposition 19.3]{ordered}: Suppose that $(P,\o)$ is naturally labeled and that every maximal chain in $P$ has length $l$. Then for all $m$, 
\begin{equation*}
\Omega(P,\o;m)=(-1)^p \Omega(P,\o; -l-m)=\Omega(P,\bar\o; l+m)
\end{equation*}
and the number of permutations in $\L(P,\o)$ with $s$ descents is equal to the number of permutations in $\L(P,\o)$ with $n-l-s$ descents.

A nice application of the reciprocity theorem for order polynomials \eqref{e-rec-ord} is Stanley's result \cite{MR0317988} that if $\chi(\lambda)$ is the chromatic polynomial of a graph $G$ with $p$ vertices then $(-1)^p\chi(-1)$ is the number of acyclic orientations of $G$. Any proper coloring of $G$ with the integers $1, 2, \dots, \lambda$ yields an acyclic orientation of $G$ in which edges are directed from the lower color to the higher. The number of proper colorings of $G$ in $\lambda$ colors associated to an  acyclic orientation $O$ is $\Omega(P_O,\o;\lambda)$ where $P$ is the poset associated to $O$ and $\o$ is a strict labeling. Then by \eqref{e-rec-ord}, $(-1)^p\Omega(P_O,\o;-1)=\Omega(P_O, \bar\o; 1) = 1$, since $\bar\o$ is a natural labeling, and thus each acyclic orientation of $G$ contributes exactly 1 to 
 $(-1)^p\chi(-1)$.

\subsection{Disjoint unions}
We may allow the labeling $\o$ of a poset $P$ to be an arbitrary function from $P$ to the positive integers as long as incomparable elements of $P$ have distinct labels. Then if $(P,\o_1)$ and $(Q,\o_2)$ are labeled posets in which the images of $\o_1$ and of $\o_2$ are disjoint, we can construct the disjoint union labeled poset $(P+Q, \o_1+\o_2)$ where the labeling $\o_1+\o_2$ is $\o_1$ on $P$ and is $\o_2$ on $Q$.  It is clear that, with the notation of Section \ref{s-theory},
\begin{equation*}
U_m(P+Q,\o_1+\o_2) = U_m(P,\o_1)U_m(Q,\o_2)
\end{equation*}
and similarly for a disjoint union of more than two posets. 
Moreover, there is a simple description of $\L(P+Q,\o_1+\o_2)$: it is the set of ``shuffles" of $\L(P,\o_1)$ and $\L(Q,\o_2)$. Thus to obtain  MacMahon's formula \eqref{e-qSN2} from \eqref{e-majdes} we take $(P,\o) = (P_1+\cdots+P_r, \o_1+\cdots +\o_r)$ where $P_i$ is a chain of 
size $p_i$ with every label equal to $i$, so that $U_m(P_i, \o_i) = \qbinom{m+p_i}{p_i}$.

Now let 
\begin{equation*}
W_s(P,\o) = \sum_{\pi} q^{\maj(\pi)},
\end{equation*}
where the sum is over all permutations $\pi\in \L(P,\o)$ with $s$ descents.
Stanley \cite[Prop.~12.6]{ordered} proves the formula
\begin{multline*}
W_s(P+Q,\o_1+\o_2) \\
= \sum_{i=0}^{|P|-1}\sum_{j=0}^{|Q|-1} q^{(s-i)(s-j)}\qbinom{|P|+j-i}{s-i}\qbinom{|Q|+i-j}{s-j}W_i(P,\o_1)W_j(Q,\o_2)
\end{multline*}
which is especially interesting in (and equivalent to) the case in which $P$ and $Q$ are chains, where it describes the enumeration by descents and major index of the shuffles of two permutations. 
Bijective proofs of this formula were later found by Goulden \cite{MR773053} and by Stadler \cite{MR1713460}.

\subsection{$\alpha$ and $\beta$}
\label{e-ab}
For any subset $S$ of $[p-1]$, let $\a(P,\o;S)$ be the number of permutations in $\L(P,\o)$ with descent set  contained in $S$ and let $\b(P,\o;S)$ be the number of permutations in $\L(P,\o)$ with descent set  equal to $S$. Then $\a(P,\o;S)=\sum_{T\subseteq S}\b(P,\o;T)$, so by inclusion-exclusion, \[\b(P,\o;S)=\sum_{T\subseteq S}(-1)^{|S|-|T|}\a(P,\o;T).\]
We can give another interpretation to $\a(P,\o;S)$. An \emph{order ideal} of $P$ is a subset $I$ of $P$ such that if $X\in I$ and $Y\prec X$ then $Y\in I$. A chain of order ideals in $P$
\begin{equation*}
\varnothing=I_0\subset I_1\subset\cdots\subset I_k=P
\end{equation*}
is called \emph{$\o$-compatible} if the restriction of $\o$ to each $I_{i+1}-I_i$ is order-preserving. (If $\o$ is natural, then any chain of order ideals is $\o$-compatible.) It is not hard to see that if the elements of $S$  are 
$m_1<m_2<\cdots<m_s$ then $\a(P,\o;S)$ is the number of $\o$-compatible chains 
$\varnothing=I_0\subset I_1\subset\cdots\subset I_{s+1}=P$ in which $|I_i|=m_i$ for $i=1,\dots, s$: given such a chain we associate to it 
the permutation consisting of the labels of $I_1$ in increasing order, followed by the labels of $I_2-I_1$ in increasing order, and so on.

We call two labelings $\o_1$ and $\o_2$ of $P$ \emph{equivalent} if $\A(P,\o_1) = \A(P,\o_2)$. Alternatively, $\o_1$ and $\o_2$ are equivalent if whenever $Y$ covers $X$ in $P$, $\o_1(X)>\o_1(Y)$ if and only if $\o_2(X)>\o_2(Y)$. If $\o_1$ and $\o_2$ are equivalent labelings, then for every $S\subseteq [p-1]$, we have $\a(P,\o_1;S)=\a(P,\o_2;S)$, and there is a simple bijection between the permutations counted by $\a(P,\o_1;S)$ and those counted by $\a(P,\o_2;S)$. It follows that $\b(P,\o_1;S)=\b(P,\o_2;S)$. This fact has interesting consequences; for example \cite[Exercise 7.95]{ec2}, it can be used  to show the existence of Solomon's descent algebra \cite{solomon} for the symmetric group.

\section{Further Developments}

\subsection{Posets}
As Stanley noted in \cite[Section 4]{ordered}, $P$-partitions are closely related to the distributive lattice $J(P)$ of order ideals of $P$, ordered by inclusion. In particular, if $\o$ is a natural labeling then the order polynomial $\Omega(P,\o;m)$ is the number of chains of order ideals
\begin{equation*}
\varnothing=I_0\subseteq I_1\subseteq\cdots\subseteq I_m=P,
\end{equation*}
and as described in Section \ref{e-ab}, $\a(P,\omega;S)$, and thus $\b(P,\omega;S)$, can be defined in terms of chains in $J(P)$. These counts of chains make sense in any graded poset with a unique minimal and maximal element, and in this context the analogue of the order polynomial is  called the \emph{zeta polynomial} and $\a$ and $\b$ are called the \emph{flag $f$-vector} and \emph{flag $h$-vector} (or \emph{rank-selected M\"obius invariant}). An account of their basic properties can be found\ in \cite[Sections 3.12 and 3.13]{ec1}. Stanley studied aspects of these concepts in  \cite{MR0309815}, \cite{MR0354472}, \cite{MR0354473}, and \cite{MR0409206}. Without further conditions the numbers $\b(S)$ need not be nonnegative, but if the edges of the Hasse diagram of the poset can be labeled with integers so that whenever $s\le t$ there is a unique saturated chain from $s$ to $t$ with nondecreasing labels (an ``R-labeling") then $\b(S)$ has a combinatorial interpretation completely analogous to the $P$-partition case. (See \cite[Section 3.14]{ec1}.)

\emph{Cohen-Macaulay posets} \cite{MR661307,MR570784} are another important class of posets for which $\b(S)$ can be shown to be nonnegative by algebraic or topological methods, but for which $\b(S)$ does not in general have a combinatorial interpretation. 

\subsection{Counting lattice points}
If $\sP$ is a lattice polytope in $\R^p$ (the convex hull of a set of lattice points) then the number of lattice points in $m\sP$ ($\sP$ dilated by a factor of $m$) is a polynomial in $m$, called the \emph{Ehrhart polynomial} of $\sP$. If $(P,\o)$ is naturally labeled then $\Omega(P,\o; m+1)$ is the Ehrhart polynomial of the \emph{order polytope} of $P$ which is the set of all $(x_1,\dots, x_p)$ in $\R^p$ satisfying $x_i\ge x_j$ whenever $i\prec j$, and $0\le x_i\le 1$ for all $i$. Some of the properties of order polynomials generalize to Ehrhart polynomials, including the reciprocity theorem for order polynomials, equation \eqref{e-rec-ord}. Stanley has made important contributions to the theory of Ehrhart polynomials and their generalizations, as described in Matthias Beck's paper \cite{beck} in this volume; see also \cite[sections 4.5--4.6]{ec1} and \cite{MR2911976}.

In \cite{MR824105}, Stanley defines the \emph{chain polytope} of the poset $P$ to be the set of points $(x_1,\dots, x_p)$ in $\R^p$ satisfying 
$x_i\ge0$ for all $i$ and $x_{i_1}+\cdots +x_{i_k}\le 1$ for every chain $i_1\prec\cdots\prec i_k$ in $P$, and proves that this polytope has the same Ehrhart polynomial as the order polytope of $P$.

\subsection{Root systems}
Gessel \cite[p.~300]{MR777705} suggested that the inequalities that define $P$-partitions could be generalized to the inequalities determined by the reflecting hyperplanes of a Coxeter group.

Victor Reiner \cite{MR1101971,reiner-signed} observed that the definition of a $P$-partition can be restated in terms of the \emph{root system} of type $A_{p-1}$, which is the set of vectors in $\R^p$ of the form $e_i-e_j$, with $i\ne j$, where $e_i$ is $i$th standard basis vector in $\R^p$. The \emph{positive roots} are the roots $e_i -e_j$ where $i<j$ and the \emph{negative roots} are the negatives of these.
Given a set $R$ of roots, we can consider the set of vectors $\s=(\s_1,\dots, \s_p)\in \N^p$ satisfying
\begin{align*}
\langle \a,\s\rangle&\ge 0, \text{ for all roots $\a$ in $R$},\\
\langle \a,\s\rangle& > 0, \text{ for all negative roots $\a$ in $R$},
\end{align*}
where $\langle\, \cdot\,, \cdot\,\rangle$ is the usual inner product on $\R^p$.
Then if this system of inequalities is consistent, the set of solutions will be the set of $(P,\o)$-partitions for some partial order on $P=[p]$ with the labeling $\o(i) = i$ for all $i\in [p]$. For example, the poset of Figure \ref{f-1} corresponds the the set of roots $\{e_2-e_1, e_2-e_3\}$; the negative root $e_2-e_1$ gives the inequality $\s_2-\s_1>0$ and the (positive) root $e_2-e_3$ gives the inequality $\s_2-\s_3\ge 0$.

Reiner \cite{MR1248063,reiner-signed} generalized this idea to arbitrary root systems, which are sets of vectors in $\R^p$ satisfying certain properties; each root system consists of a set of positive roots and their negatives, the negative roots. To each root system is associated its \emph{Weyl group}, which is the finite group of isometries of $\R^p$ generated by the reflections in the roots. (So for the root system of type $A_{p-1}$, the Weyl group is the symmetric group ${\mathfrak S}_p$.) Reiner shows that there is a generalization of the fundamental theorem of $P$-partitions to root systems, in which the role of the symmetric group is replaced by the corresponding Weyl group.

He then studies in particular the case of the root system of type $B_p$, in which the positive roots are $e_i$ for $1\le i\le p$ and $e_i+e_j$ and $e_i-e_j$ for $1\le i<j\le n$. Thus, for example, if we take the roots $e_1$, $-e_2-e_3$, and $e_3-e_1$, then the corresponding inequalities are $\s_1\ge0$, $-\s_2-\s_3>0$ and $\s_3 -\s_1<0$. (Unlike the case of ordinary $P$-partitions, here we allow the $\s_i$ to take on arbitrary integer values.) The elements of the Weyl group of type $B_p$, the hyperoctahedral group, may be viewed as \emph{signed permutations}, which are permutations $\pi$ of the set
$\pm[p]=\{-p,\cdots, -1, 1, \cdots p\}$,  such that $\pi(-i) = -\pi(i)$ for all $i\in \pm[p]$; a signed permutation is determined by its values on $[p]$. The general definition of descent set for a Weyl group reduces in this case to 
\begin{equation*}
\S(\pi) = \{\,i \in [p] \mid \ \pi(i)>\pi(i+1)\,\},
\end{equation*}
where we take $\pi(p+1) =p+1$ and use the order\footnote{Different choices for the roots would allow similar results with the usual order on $\pm[p]$.} 
$1<2<\cdots < p+1 <-p< \cdots <-2<-1$. Reiner then obtains for signed permutations analogues of all the basic $P$-partition for ordinary permutations.

Chak-On Chow \cite[Section 2]{MR2717011} studied ``$P$-partitions of type $B$" with a closely related, but somewhat different approach,  using ``type $B$ posets" which are partial orders $\prec$ on the set $\{-p,\cdots, -1, 0, 1,\cdots,p\}$ with the property that $i\prec j$ if and only if $-j \prec -i$.

Further work on root system analogues of $P$-partitions was undertaken by John Stembridge in his study of Coxeter cones  \cite{MR2388082}.

\subsection{Lexicographic inequalities}
MacMahon, in \cite{multi} and other papers (see \cite[Chapter 8]{MR514405}), studied \emph{multipartite partitions} which are expressions of ``multipartite numbers" as sums of multipartite numbers, where a multipartite number is a tuple of nonnegative integers. The number of partitions of the multipartite number $(n_1,\dots, n_s)$ with $p$ parts, where $(0,\dots,0)$ is allowed as a part, is the coefficient of $q_1^{n_1}\cdots q_{s}^{n_s}$ in $\phi_r(x)$, where 
\begin{equation*}
\sum_{p=0}^\infty \phi_p(q_1,\dots, q_s)z^p = \prod_{k_1,\dots, k_s=0}^\infty \frac{1}{1-q_1^{k_1}\cdots q_s^{k_s}z}
\end{equation*}
It is not hard to show that there exist polynomials $\Lambda_p(q_1,\dots, q_s)$ such that 
\begin{equation*}
\phi_p(q_1,\dots, q_s) = \frac{\Lambda_p(q_1,\dots, q_s)}{(q_1;q_1)_p\cdots (q_s;q_s)_p}
\end{equation*}
where $(q;q)_s=(1-q)\cdots (1-q^s)$. E. M. Wright \cite{MR0084012} conjectured that the polynomials $\Lambda_p(q_1,\dots, q_s)$ have nonnegative coefficients, and his conjecture was proved by Basil Gordon \cite{MR0157959} in 1963. We illustrate Gordon's approach with the simplest example, when $p=s=2$. An unordered pair $(a_1,b_1), (a_2,b_2)$ of bipartite numbers may be arranged in decreasing lexicographic order, so counting such pairs is equivalent to counting solutions of the lexicographic inequality $(a_1,b_1)\ge (a_2, b_2)$. The set of solutions of this lexicographic inequality is the disjoint union of the solutions of 
\begin{gather*}
a_1\ge a_2, \ b_1\ge b_2\\
\shortintertext{and}
a_1>a_2, \ b_1 < b_2.\\
\end{gather*}
The solutions of the first inequalities contribute $1/(q_1;q_1)_2(q_2;q_2)_2$ to $\Lambda_2(q_1,q_2)$ and the solutions of the second inequalities contribute $q_1q_2/(q_1;q_1)_2(q_2;q_2)_2$, so $\Lambda_2(q_1,q_2) = 1+q_1q_2$. 

Gordon proved Wright's conjecture by showing  that the lexicographic inequalities specifying the terms in $\phi_p(q_1,\dots, q_s)$ can in general be decomposed in a similar way. D. P. Roselle \cite{MR0342406}, using essentially the same approach,  gave a simple combinatorial  interpretation to the coefficient of $q_1^{i_1}q_2^{i_2}$ in $\Lambda_r(q_1,q_2)$; it is the number of permutations $\pi$ of $[p]$ such that $\maj(\pi) = i_1$ and $\maj(\pi^{-1})=i_2$. Garsia and Gessel \cite{MR532836} showed that, more generally, the coefficient of $q_1^{i_1}\cdots q_s^{i_s}$ in $\Lambda_p(q_1,\cdots,q_s)$ is the number of $s$-tuples $(\pi_1,\dots, \pi_s)$ of permutations of $[p]$ whose product $\pi_1\cdots \pi_s$ is the identity permutation such that $\maj(\pi_j)=i_j$ for each $j$. They also showed that by considering multipartite partitions with bounded part sizes, one can count these $s$-tuples of permutations by major index and number of descents. Further results along these lines have been found by a number of authors \cite{MR2780854,  MR2195428, MR2196519, moynihan,MR2181371,MR1101971,MR1248063}. 

Gessel \cite{MR777705} studied ``multipartite $P$-partitions" in which the inequalities defining $(P,\o)$-partitions are applied to multipartite numbers ordered lexicographically; his results enumerate $s$-tuples of permutations whose product is in $\L(P,\o)$ by their descent sets.

Another application of lexicographic inequalities related to $P$-partitions was given by Gessel and Reu\-ten\-auer \cite{MR1245159}. A \emph{Lyndon word} is a sequence of nonnegative integers that is lexicographically strictly less than all of its cyclic permutations. Thus  the word $(a_1, a_2, a_3)$ is a Lyndon word if and only if $(a_1, a_2, a_3) < (a_2, a_3, a_1)$ and $(a_1, a_2, a_3)<(a_3, a_1, a_2)$. The set of solutions of these lexicographic inequalities is the disjoint union of the solutions of 
\begin{gather*}
a_1\le a_2 < a_3\\
\shortintertext{and}
a_1 < a_3 \le a_2.
\end{gather*}
Gessel and Reutenauer showed that a similar decomposition exists for Lyndon words of any length, and more generally, for multisets of Lyndon words, and this allowed them to count permutations of a given cycle type by their descent sets. A generalization to hyperoctahedral groups was given by St\'ephane Poirier \cite{MR1603753}.

\subsection{Quasi-symmetric functions}
\label{s-qs}
In his memoir \cite{ordered}, Stanley considered the generating function for $(P,\o)$-partitions
\begin{equation*}
F(P,\o)=\sum_{\s\in \A(P,\o)}x_1^{\s(1)}x_2^{\s(2)}\cdots x_p^{\s(p)}
\end{equation*}
in which different $(P,\o)$-partitions contribute different terms.  The theory of Schur functions (see, for example, \cite[Section 7.10]{ec2}) suggests looking at the less refined generating function 
\begin{equation}
\label{e-qsgf}
\Gamma(P,\o) = \sum_{\s\in \A(P,\o)}x_{\s(1)}x_{\s(2)}\cdots x_{\s(p)}
\end{equation}
discussed briefly by Stanley \cite[p.~81]{ordered} and  in more detail by 
Gessel\footnote{Gessel took $P$-partitions to be order-preserving, rather than order-reversing, so his $\Gamma(P,\o)$ is slightly different from that defined here. } \cite{MR777705}. (See also \cite[Section 7.19]{ec2}.)
By the fundamental theorem of $P$-partitions, 
\begin{equation*}
\Gamma(P,\o)  = \sum_{\pi\in \L(P,\o)} \Gamma(\pi,\o).
\end{equation*}
The \emph{fundamental quasi-symmetric functions}, denoted $F_\a$ or $L_\a$, are indexed by compositions (sequences of positive integers) and defined 
as follows: if $\a=(\a_1,\dots a_k)$ is a composition of $p$ then 
\begin{equation*}
F_\a = \sum x_{i_1} x_{i_2}\cdots x_{i_p}
\end{equation*}
where the sum is over all $i_1\le i_2\le\cdots\le i_p$ satisfying $i_j < i_{j+1}$ if $j\in \{\a_1, \a_1+\a_2,\dots, \a_1+\cdots+\a_{k-1}\}$. It is not hard to show that the $F_\a$ are linearly independent and generate a ring, called the \emph{ring of quasi-symmetric functions}, that contains the ring of symmetric functions \cite[Chapter 7]{ec2}. Thus the information contained in $\Gamma(P,\o)$ is precisely the multiset of descent sets of the 
permutations in $\L(P,\o)$; i.e., numbers $\b(P,\omega;S)$. The quasi-symmetric generating function \eqref{e-qsgf} extends to an encoding of the flag $h$-vector of a graded poset, as studied by Richard Ehrenborg \cite{MR1383883}.

The theory of quasi-symmetric functions has proven useful in a number of enumeration problems. For example, Stanley \cite{MR782057} used them to define what are now called ``Stanley symmetric functions" in the study of reduced  decompositions in symmetric groups.

There is a comultiplication on the ring of quasi-symmetric functions that makes it into a Hopf algebra (and an ``internal" comultiplication that makes it a bialgebra). The dual Hopf algebra is the algebra of \emph{non-commutative symmetric functions} that has been studied extensively by Jean-Yves Thibon and others in a series of papers beginning with 
\cite{MR1327096}; see also Malvenuto and Reutenauer \cite{MR1358493}. Type $B$ quasi-symmetric functions and noncommutative symmetric functions were studied by Chow \cite{MR2717011}. Another related algebra with enumerative applications is the Malvenuto-Reutenauer algebra \cite{MR2103213,malvenuto,MR1358493}.

A different encoding of the flag $h$-vector of a poset is the \emph{$\mathbf{ab}$-index}, which is especially useful in studying Eulerian posets \cite{MR1651249,MR1283084}.

\subsection{Enriched $P$-partitions}
John Stembridge \cite{enriched} introduced a generalization of $(P,\o)$-partitions that interpolates between  $(P,\o)$-partitions and $(P,\bar\o)$-partitions. We introduce the following total ordering on the set $\P'$ of nonzero integers:
\begin{equation*}
-1< +1 < -2 < +2 < -3 < +3<\cdots;
\end{equation*}
for $k\in \P'$ the notations $k>0$ and $|k|$ retain their usual meanings.
Then an enriched $(P,\o)$-partition is a map $\s: P\to \P'$ such that
 for all $X\prec Y$ in $P$ we have%
\footnote{Stembridge defined enriched $P$-partitions to be order-preserving; for consistency we define them here to be order-reversing.}
\begin{enumerate}
\item[(i)]  $\s(X)\ge\s(Y)$
\item[(ii)] If $\s(X) =\s(Y)>0$ then $\o(X) <\o (Y)$ 
\item[(iii)] If $\s(X) = \s(Y) <0$ then $\o(X) > \o(Y)$.
\end{enumerate}
Note that if the image of $\s$ lies in $\{+1,+2,\dots\}$ then (iii) is vacuous and (ii) is equivalent to the condition that if $\o(X)>\o(Y)$ then $\s(X)>\s(Y)$, so $\s$ is an ordinary $(P,\o)$-partition, and if the image of $\s$ lies in $\{-1, -2, \cdots \}$ then (ii) is vacuous and (iii) is equivalent to the condition that if $\o(X)<\o(Y)$ then $\s(X) > \s(Y)$, so $\s$ is essentially a $(P,\bar\o)$-partition. Stembridge proves a version of the fundamental theorem for enriched $(P,\o)$-partitions, and defines the quasi-symmetric generating function 
\begin{equation*}
\Delta(P,\o) = \sum_\s \prod_{X\in P}x_{|\s(X)|},
\end{equation*}
so by the fundamental theorem, 
\begin{equation*}
\Delta(P,\o) = \sum_{\pi\in \L(P,\o)} \Delta(\pi,\o).
\end{equation*}
It is a remarkable fact that $\Delta(\pi,\o)$ depends only on the \emph{peak set} of $\pi$, that is, the set $\{\, i \mid \pi(i-1) < \pi(i) > \pi(i+1)\,\}$. The distinct $\Delta(\pi,\o)$ form a basis for a subalgebra of the algebra of quasi-symmetric functions.

Stembridge also discusses the analogue of the order polynomial for enriched $P$-partitions and studies cases in which $\Delta(P,\o)$ is symmetric, which are related to Schur's $Q$-functions. 

Kathryn Nyman \cite{MR2001673} used enriched $P$-partitions to prove the existence of the ``peak algebra" of the symmetric group.  T.~Kyle Petersen \cite{MR2296309} studied type $B$ enriched $P$-partitions and applied them to type $B$ peak algebras.
Enriched $P$-partitions have also been applied to the study of chains in Eulerian posets \cite{MR1982883}.

\subsection{Additional applications and developments}

SeungKyung Park \cite{park} studied naturally labeled posets $P$ whose order polynomial $\Omega_P(n)$ is the Stirling number of the second kind $S(k+n,n)$, thereby giving a new combinatorial interpretation and $q$-analogue to the polynomials $B_k(t)$ defined by 
\begin{equation*}
\sum_{n=0}^\infty S(k+n,n) t^n = \frac{B_k(t)}{(1-t)^{2k+1}}.
\end{equation*}
Combinatorial interpretations for these polynomials had been given earlier by John Riordan \cite{MR0429582} and by Gessel and Stanley \cite{MR0462961}.

Sangwook Ree \cite{sree} applied $P$-partitions to count restricted lattice paths in the plane by left turns, obtaining generalizations of $q$-Narayana numbers. 
Petter Br\"and\'en \cite{MR2047757}, taking a similar approach, gave several interpretations to $q$-Narayana numbers in counting Dyck paths.

Joseph Neggers \cite{MR0551484} conjectured in 1978 that for any naturally labeled poset $(P,\o)$ the polynomial
$\sum_{\pi\in \L(P,\o)} t^{\des(\pi)}$ has all real roots. In 1986, Stanley conjectured that this holds more generally for any labeled poset $(P,\o)$ (see 
\cite{MR963833}). Stanley's conjecture was disproved by Petter Br\"and\'en \cite{MR2119757} in 2004 and Neggers's conjecture was disproved by John Stembridge 
\cite{MR2262844}
in 2007.

Peter McNamara and Christophe Reutenauer \cite{MR2195427} used $P$-partitions to study idempotents in the group algebra of the symmetric group. 

McNamara and  Ryan Ward \cite{MR3245895} studied the question of when two different labeled posets have the same generating function \eqref{e-qsgf}.

Valentin F\'eray and Victor Reiner \cite{MR2913529} explored  connections between $P$-par\-ti\-tions and commutative algebra, and in particular described a class of naturally labeled  posets for which the sum $\sum_{\pi\in \L(P)}q^{\maj(\pi)}$ factors nicely.  

Lo\"\i c Foissy and Claudia Malvenuto \cite{malv} reinterpreted the fundamental theorem of $P$-partitions as an injective Hopf algebra morphism and generalized it to pre-orders, leading to a Hopf algebra on finite topologies.

\subsection*{Acknowledgments} I would like to thank
Christian Krattenthaler,
Claudia Malvenuto, 
T. Kyle Petersen, 
Victor Reiner, and
John Stembridge
for their helpful suggestions, and Richard Stanley for his development of the theory of $P$-partitions.

\providecommand{\bysame}{\leavevmode\hbox to3em{\hrulefill}\thinspace}

\end{document}